\documentclass{amsart}[12pt]

\usepackage{amscd,amssymb}
\usepackage{comment}

\newtheorem{theorem}{Theorem}[section]
\newtheorem{lemma}[theorem]{Lemma}

\theoremstyle{definition}
\newtheorem{definition}[theorem]{Definition}
\newtheorem{example}[theorem]{Example}

\newtheorem{remark}[theorem]{Remark}

\newcommand{\GL} {\mathrm{GL}}
\newcommand{\SL} {\mathrm{SL}}
\newcommand{\OO} {\mathrm{O}}
\newcommand{\SO} {\mathrm{SO}}
\newcommand{\Sp} {\mathrm{Sp}}

\def\CC {{\mathbb C}}     
\def\NN {{\mathbb N}}     
\def\ZZ {{\mathbb Z}}     

\begin{document}

\title[Hilbert series and invariants in exterior algebras]
{Hilbert series and invariants in exterior algebras}

\author[Elitza Hristova]
{Elitza Hristova}
\address{Institute of Mathematics and Informatics,
Bulgarian Academy of Sciences,
Acad. G. Bonchev Str., Block 8,
1113 Sofia, Bulgaria}
\email{e.hristova@math.bas.bg}


\subjclass[2010] {13A50; 15A72; 15A75; 20G05.}

\keywords{invariant theory, Hilbert series, multiplicity series, exterior algebras}

\begin{abstract}In this paper, we consider the exterior algebra $\Lambda(W)$ of a polynomial $\GL(n)$-module $W$ and use previously developed methods to determine the Hilbert series of the algebra of invariants $\Lambda(W)^G$, where $G$ is one of the classical complex subgroups of $\GL(n)$, namely $\SL(n)$, $\OO(n)$, $\SO(n)$, or $\Sp(2d)$ (for $n=2d$). Since $\Lambda(W)^G$ is finite dimensional, we apply the described method to compute a lot of explicit examples. For $\Lambda(S^3\CC^3)^{\SL(3)}$, using the computed Hilbert series, we obtain an explicit set of generators.
\end{abstract}

\maketitle

\section{Introduction}
Let $\GL(n)$ denote the complex general linear group and let $W$ be a finite dimensional polynomial $\GL(n)$-module. Then $W$ can be written as a direct sum of its irreducible components
\begin{equation*} 
W = \bigoplus_{\lambda} k(\lambda) V^n_{\lambda},
\end{equation*}
where $\lambda = (\lambda_1, \ldots, \lambda_n) \in \NN_0^n$, $\lambda_1\geq \lambda_2 \geq \cdots \geq \lambda_n \geq 0$,
is a non-negative integer partition (here $\NN_0=\{0,1,2,\ldots\}$) and
$V^n_{\lambda}$ denotes the irreducible $\GL(n)$-module with highest weight $\lambda$.
Let $A = \bigoplus_{i \geq 0} A^i$ be a finitely generated graded algebra over $\CC$ for which each homogeneous component $A^i$ is a polynomial $\GL(n)$-module. In the papers \cite{DH1, DH2}, Drensky and the author develop a method for determining the Hilbert series of the algebra of invariants $A^G$, where $G$ is one of the classical complex groups $\OO(n)$, $\SO(n)$, and $\Sp(2d)$ (the last in the case when $n = 2d$).
These results can be easily extended to the case $G = \SL(n)$ using an earlier work of Drensky and other authors \cite{BBDGK}.
A general class of examples of an algebra $A$ with the above properties is given by the quotient $T(W)/I$, where $T(W)$ denotes the tensor algebra of a polynomial $\GL(n)$-module $W$ and $I$ is a $\GL(n)$-invariant ideal in $T(W)$. This class includes the symmetric and exterior algebras of $W$, $S(W)$ and $\Lambda(W)$, as well as the class of relatively free algebras in varieties of associative algebras. In the papers \cite{DH1, DH2}, Drensky and the author 
use the developed method to study the Hilbert series $H(A^G,t)$ for $A$ being the symmetric algebra $S(W)$ and for certain relatively free algebras 
with the property that their varieties are the only minimal varieties of exponent $2$. In the present paper, we consider the case $A = \Lambda(W)$ and study the algebra of invariants $\Lambda(W)^G$, where $G$ is again one of $\SL(n)$, $\OO(n)$, $\SO(n)$, and $\Sp(2d)$. 
The case $A = \Lambda(W)$ is remarkable with the property that the considered Hilbert series are polynomials, which allows for explicit computation of a lot of examples. 
Invariants in exterior algebras for different actions of classical groups and related problems are the topic of many recent papers, see, e.g., \cite{Do1, Do2, I}.

The present paper is organized as follows. In Section \ref{sec_prelim} we give a brief introduction to the general method for determining Hilbert series which is developed in \cite{DH2}. In Section \ref{sec_method} we describe how the general method is applied to $A=\Lambda(W)$ and in Section \ref{sec_examples} we take concrete examples of $W$ and compute explicitly the respective Hilbert series. The results are given in Tables 1-9 and in Example \ref{ex_Hilb5}. In some cases, we use the computed Hilbert series to describe also a set of generators of the respective algebra of invariants (see Example \ref{ex_LambdaS3_SL_inv}). Since for most of the examples we use computer programs written in Mathematica, in the Appendix we provide a part of the source code.

\section{Preliminaries} \label{sec_prelim}
The ground field is $\CC$. We start by recalling the notion of Hilbert series, which applies to finitely generated graded algebras (resp., vector spaces) over $\CC$.
\begin{definition} \label{def_HilbSeries}
Let $\displaystyle A = \bigoplus_{i\geq 0} A^i$ be a finitely generated graded algebra (or a vector space) over $\CC$
such that $A^0=\CC$ or $A^0=0$.
The {\bf Hilbert series} of $A$ is the formal power series
\[
H(A, t) = \sum_{i\geq 0} (\dim A^i) t^i.
\]
\end{definition}

The next definition is a generalization of the notion of Hilbert series for multigraded algebras (resp., for multigraded vector spaces).
\begin{definition} \label{def_HilbSeriesMult} Let
\[
A = \bigoplus_{\mu \in \NN_0^n} A(\mu)
\]
be a finitely generated algebra or vector space with an $\NN_0^n$-grading. The Hilbert series of $A$ with respect to this grading is the formal power series $H(A, x_1, \ldots, x_n) \in \ZZ[[x_1, \dots, x_n]]$ defined by
\[
H(A, x_1, \ldots, x_n) = \sum_{\mu = (\mu_1, \ldots, \mu_n) \in \NN_0^n} \dim A(\mu) x_1^{\mu_1} \cdots x_n^{\mu_n}.
\]
\end{definition}

One example of a vector space with an $\NN_0^n$-grading is the $\GL(n)$-module $V^n_{\lambda}$ together with its weight space decomposition. The Hilbert series of $V^n_{\lambda}$ with respect to this grading coincides with the character $\chi_{V^n_{\lambda}}(x_1, \dots, x_n)$ of $V^n_{\lambda}$ and consequently has the form
\[
H(V^n_{\lambda}, x_1, \ldots, x_n) = \chi_{V^n_{\lambda}}(x_1, \dots, x_n) = S_{\lambda}(x_1, \ldots, x_n),
\]
where $S_{\lambda}(x_1, \ldots, x_n)$ is the Schur polynomial corresponding to the partition $\lambda$. Similarly, each polynomial $\GL(n)$-module $W$ also has an $\NN_0^n$-grading and a corresponding Hilbert series which is equal to the character of $W$.

Next, as in the Introduction, we consider the special class of finitely generated graded algebras $\displaystyle A = \bigoplus_{i\geq 0} A^i$, such that each homogeneous component $A^i$ for $i \geq 0$ is a polynomial $\GL(n)$-module. Equivalently, $A$ can be defined as an arbitrary direct sum of polynomial $\GL(n)$-modules. Then, we can write the decomposition of $A$ in the following way:
\[
A = \bigoplus_{i\geq 0} A^i = \bigoplus_{i \geq 0} \bigoplus_{\lambda} m_i(\lambda)V_{\lambda},
\]
where $A^i = \bigoplus_{\lambda} m_i(\lambda)V_{\lambda}$ is the decomposition of $A^i$ into irreducible $\GL(n)$-modules and the sum runs over all partitions $\lambda$ in $\NN_0^n$. Hence, $A$ has two gradings -- the $\NN_0$-grading coming from the decomposition into homogeneous components and an $\NN_0^n$-grading coming from the $\NN_0^n$-grading of each homogneous component. 
Following \cite{BBDGK, DH1, DH2}, we introduce a Hilbert series of $A$ which takes into account both gradings:
\[
H(A; x_1, \dots, x_n, t) = \sum_{i \geq 0} H(A^i, x_1, \ldots, x_n) t^i = \sum_{i \geq 0}\left (\sum_{\lambda} m_i(\lambda) S_{\lambda}(x_1, \dots, x_n)\right ) t^i.
\]
Let $\mathrm{S}_n$ denote the symmetric group in $n$ variables. Then, $H(A; x_1, \dots, x_n, t) \in \ZZ[[x_1, \dots, x_n]]^{\mathrm{S}_n}[[t]]$, i.e., the coefficient in front of $t^i$ for each $i \geq 0$ is a symmetric function in the variables $x_1, \dots, x_n$. For any symmetric function $f(x_1, \dots, x_n) \in \CC[[x_1, \dots, x_n]]^{\mathrm{S}_n}$, one can define the notion of multiplicity series (see, e.g., \cite{BBDGK}).
Following \cite{BBDGK, DH1, DH2}, we introduce a generalization of the notion of multiplicity series for functions $f(x_1, \dots, x_n, t)$ in $\CC[[x_1, \dots, x_n]]^{\mathrm{S}_n}[[t]]$. Let $f(x_1, \dots, x_n, t) \in \CC[[x_1, \dots, x_n]]^{\mathrm{S}_n}[[t]]$, i.e., $f(x_1, \dots, x_n, t)$ can be written as
\[ 
f(x_1, \dots, x_n, t) = \sum_{i \geq 0} \left(\sum_{\lambda} c_i(\lambda)S_{\lambda}(x_1, \dots, x_n)\right)t^i, 
\]
where $c_i(\lambda) \in \CC$ and the second sum runs over all partitions $\lambda \in \NN_0^n$. The multiplicity series of $f(x_1, \dots, x_n, t)$ is defined by:
\begin{equation} \label{eq_multSeries1}
M(f; x_1, \dots, x_n, t) = \sum_{i \geq 0}\left(\sum_{\lambda} c_i(\lambda)x_1^{\lambda_1} \cdots x_n^{\lambda_n} \right)t^i.
\end{equation}
By the change of variables $v_1 = x_1$, $v_2 = x_1x_2$, \dots, $v_n = x_1\cdots x_n$, we rewrite the above multiplicity series in the following form:
\begin{align} \label{eq_multSeries2}
M'(f; v_1,\ldots, v_n, t) = \sum_{i \geq 0}\left(\sum_{\lambda}
c_i(\lambda)v_1^{\lambda_1 - \lambda_2} v_2^{\lambda_2 - \lambda_3} \cdots v_{n-1}^{\lambda_{n-1} - \lambda_n} v_n^{\lambda_n}\right)t^i.
\end{align}

For any subgroup $G$ of $\GL(n)$, the algebra of invariants $A^G$ inherits the $\NN_0$-grading of $A$ and therefore one can study the Hilbert series $H(A^G, t)$.
The following theorem, given in \cite{DH2}, shows the relation between $H(A^G,t)$ and the multiplicity series of $H(A; x_1, \dots, x_n, t)$ defined by (\ref{eq_multSeries1}) and (\ref{eq_multSeries2}).

\begin{theorem}[\cite{DH2}] \label{thm_HilbertSeries} Let $\displaystyle A = \bigoplus_{i\geq 0} A^i$ be a finitely generated graded algebra, such that each homogeneous component $A^i$ for $i \geq 0$ is a polynomial $\GL(n)$-module.

{\rm (i)} The Hilbert series of the algebra of invariants $A^{\SL(n)}$ is given by
\[
H(A^{\SL(n)}, t) = M'(H(A); 0,\ldots, 0, 1, t).
\]


{\rm (ii)} The Hilbert series of $A^{\Sp(2d)}$ (where $n = 2d$) is given by
\[
H(A^{\Sp(2d)}, t) = M'(H(A); 0, 1, 0, 1 ,\ldots, 0, 1, t).
\]

{\rm (iii)} The Hilbert series of $A^{\OO(n)}$ is
\[
H(A^{\OO(n)}, t) = M_{n}(t),
\]
where $M_n$ is defined iteratively in the following way:
\[
M_1(x_2, \ldots, x_n, t) = \frac{1}{2}\left(M(H(A); -1, x_2, \ldots, x_n,t) + M(H(A); 1, x_2, \ldots, x_n, t)\right),
\]
\[
M_2(x_3, \ldots, x_n, t) = \frac{1}{2}(M_1(-1, x_3, \ldots, x_n, t) + M_1(1, x_3, \ldots, x_n, t))
\]
\[
\cdots\cdots\cdots
\]
\[
M_n(t) = \frac{1}{2}(M_{n-1}(-1, t) + M_{n-1}(1, t)).
\]

{\rm (iv)} The Hilbert series of $A^{\SO(n)}$ is
\[
H(A^{\SO(n)}, t) = M'_{n}(t),
\]
where
\[
M'_1(v_2, \ldots, v_n, t) = \frac{1}{2}(M'(H(A); -1, v_2, \ldots, v_n,t) + M'(H(A); 1, v_2, \ldots, v_n, t)),
\]
\[
M'_2(v_3, \ldots, v_n, t) = \frac{1}{2}(M'_1(-1, v_3, \ldots, v_n, t) + M'_1(1, v_3, \ldots, v_n, t))
\]
\[
\cdots\cdots\cdots
\]
\[
M'_{n-1}(v_n, t) = \frac{1}{2}(M'_{n-2}(-1,v_n, t) + M'_1(1, v_n, t)),
\]
\[
M'_n(t) = M'_{n-1}(1, t).
\]
\end{theorem}

\section {The algebra of invariants $\Lambda(W)^G$} \label{sec_method}
In this section we show how Theorem \ref{thm_HilbertSeries} can be applied to determine the Hilbert series of $A^G$ for $A = \Lambda(W)$.

Let $W$ be a $p$-dimensional polynomial $\GL(n)$-module, let $\xi_1, \dots, \xi_p$ be a basis of $W$ consisting of weight vectors, and let $\alpha_1 =(\alpha_{11}, \dots, \alpha_{1n})$, $\dots$, $\alpha_p =(\alpha_{p1}, \dots, \alpha_{pn})$ denote the corresponding weights (with possible repetitions). The Hilbert series of $W$ with respect to the $\NN_0^n$-grading given by the weight space decomposition of $W$ is equal to the character $\chi_W$ of $W$ and hence
\[
H(W, x_1, \dots, x_n) = \sum_{j =1}^p x_1^{\alpha_{j1}} \cdots x_n^{\alpha_{jn}}.
\]

Similarly, the Hilbert series $H(\Lambda^k(W), x_1, \dots, x_n)$ of $\Lambda^k(W)$ for $0\leq k \leq p$ is equal to the character $\chi_{\Lambda^k(W)}$. It is a standard fact (see, e.g., \cite{FH}) that if $W = V$, where $V$ denotes the natural representation of $\GL(n)$, the character of $\Lambda^k(V)$ is equal to the $k$-th elementary symmetric polynomial, i.e.,
\[
\chi_{\Lambda^k(V)}(x_1, \dots, x_n) =  E_k(x_1, \dots, x_n) = \sum_{1\leq s_1 < \dots< s_k  \leq n} x_{s_1} \cdots x_{s_k}.
\]
Using the above formula, one shows that 
\[
H(\Lambda^k(W), x_1, \dots, x_n) = \sum_{1\leq s_1 < \dots< s_k \leq n}(x_1^{\alpha_{s_11}}\cdots x_n^{\alpha_{s_1n}})\cdots (x_1^{\alpha_{s_k1}}\cdots x_n^{\alpha_{s_kn}}).
\]
Therefore,
\begin{align} \label{eq_HilbSeriesLambda2}
H(\Lambda(W), x_1, \dots, x_n,t) = \prod_{j=1}^p(1+x_1^{\alpha_{j1}}\cdots x_n^{\alpha_{jn}}t).
\end{align}

Lemma 3 from \cite{Be} describes a way to determine the multiplicity series of each symmetric function $f(x_1, \dots, x_n) \in \CC[[x_1, \dots, x_n]]^{\mathrm{S}_n}$ (see also Lemma 1.1 from \cite{BBDGK}). The following lemma is a direct generalization for functions in $\CC[[x_1, \dots, x_n]]^{\mathrm{S}_n}[[t]]$.

\begin{lemma}\label{lemma_MultSeries}
Let $X = \{x_1, \dots, x_n\}$ and let $f(X,t) \in \CC[[X]]^{\mathrm{S}_n}[[t]]$. Let
\[
g(X,t) = f(X,t)\prod_{i<j}(x_i - x_j) = \sum_{i \geq 0} (\sum_{r_{i_j} \geq 0} \alpha_i(r_{i_1}, \dots, r_{i_n})x_1^{r_{i_1}} \cdots x_n^{r_{i_n}})t^i,
\]
for some $\alpha_i(r_{i_1}, \dots, r_{i_n}) \in \CC$. Then the multiplicity series of $f(X,t)$ is given by
\begin{align*}
&M(f; x_1, \dots, x_n,t) =\\
&\frac{1}{x_1^{n-1}x_2^{n-2} \cdots x_{n-2}^2 x_{n-1}}\sum_{i \geq 0}( \sum_{r_{i_j} > r_{i_{j+1}}} \alpha_i(r_{i_1}, \dots, r_{i_n})x_1^{r_{i_1}} \cdots x_n^{r_{i_n}})t^i,
\end{align*}
where the sum is over all $r_i = (r_{i_1},\dots, r_{i_n})$ such that $r_{i_1} > r_{i_2} > \cdots > r_{i_n}$.
\end{lemma}

In the next section, we use this general procedure for concrete choices of $W$.

\section{Examples} \label{sec_examples}
Let $V = \CC^n$ denote the natural $\GL(n)$-module and let again $G$ be one of $\SL(n)$, $\OO(n)$, $\SO(n)$, and $\Sp(2d)$ (for $n = 2d$). In the first set of examples, we consider the module $W = S^kV$, for $k\geq 3$, and determine the Hilbert series $H(\Lambda(W)^G, t)$ for fixed values of $n$ and $k$. The case $k=2$ is not included in the consideration, since explicit expressions for $H(\Lambda(S^2V)^G, t)$ and $H(\Lambda(\Lambda^2V)^G, t)$ for all $n$ are already obtained in \cite{DH1}. 

We notice first that when $W = S^kV$, Equation (\ref{eq_HilbSeriesLambda2}) implies that
\[
H(\Lambda(S^kV), x_1, \dots, x_n,t) = \prod_{i_1 + \dots + i_n = k} (1 + x_1^{i_1}\cdots x_n^{i_n}t).
\]

In the simplest case, when $n=2$ and $k=3$, we show the computation of $H(\Lambda(S^kV)^G,t)$ in the following example. 
\begin{example}
Let $n = 2$ and $W = S^3V$, where $V = \CC^2$. 
Then
\[
H(\Lambda(S^3V), x_1, x_2,t) = (1+x_1^3t)(1+x_2^3t)(1+x_1^2x_2t)(1+x_1x_2^2t).
\]
Using Lemma \ref{lemma_MultSeries} we obtain that
\[
M(H(\Lambda(S^3V)); x_1, x_2,t) = 1+ x_1^3t + x_1^3x_2^3t^2 + x_1^5x_2t^2 + x_1^6x_2^3t^3  + x_1^6x_2^6t^4.
\]
By definition, the multiplicity series $M(H(\Lambda(S^3V)); x_1, x_2,t)$ carries the information about the $\GL(2)$-structure of $\Lambda(S^3V)$. Hence, in the notations from the Introduction, we obtain the following decomposition of $\Lambda(S^3V)$ as a $\GL(2)$-module
\[
\Lambda(S^3V) \cong \CC \oplus S^3V \oplus (V^2_{(3,3)} \oplus V^2_{(5,1)}) \oplus V^2_{(6,3)} \oplus V^2_{(6,6)}.
\]
In particular,
$\Lambda^2(S^3V)\cong V^2_{(3,3)} \oplus V^2_{(5,1)}$, $\Lambda^3(S^3V)\cong V^2_{(6,3)}$ and $\Lambda^4(S^3V)\cong V^2_{(6,6)}$.

Next, using Theorem \ref{thm_HilbertSeries} we determine the following expressions for $H(\Lambda(S^3V)^G, t)$:
\[ 
H(\Lambda(S^3V)^{\SL(2)}, t) = H(\Lambda(S^3V)^{\Sp(2)}, t) = 1+ t^2 + t^4;
\]
\[ 
H(\Lambda(S^3V)^{\OO(2)}, t) = 1+ t^4;
\]
\[ 
H(\Lambda(S^3V)^{\SO(2)}, t) = 1+ 2t^2 + t^4.
\] 
\end{example}

For computing further examples we have written short programs in Mathematica (respectively for $n=2$,$3$,$4$, and $5$ variables) and in Tables 1-6 below and in Example \ref{ex_Hilb5} we give the corresponding expressions for the Hilbert series. 

\begin{table}[ht] \label{table1}
\caption{Hilbert series for $n=2$}
\begin{tabular} {|c|l|}
\hline
$k$ & $H(\Lambda(S^kV)^{\SL(2)}, t)$ \\
\hline
$3$ & $1+ t^2 + t^4$\\
\hline
$4$ & $1+ t^5$\\
\hline
$5$ & $1+t^2+t^4+t^6$\\
\hline
$6$ & $1+t^3+t^4+t^7$\\
\hline
$7$ & $1+t^2+t^4+t^6+t^8$\\
\hline
$8$ & $1+t^4+t^5+t^9$\\
\hline
$9$ & $1+t^2+2t^4+2t^6+t^8+t^{10}$\\
\hline
$10$ & $1+t^3+t^4+t^7+t^8+t^{11}$\\
\hline
$11$ & $1+t^2+2t^4+3t^6+2t^8+t^{10}+t^{12}$\\
\hline
$12$ & $1+2t^4+2t^5+2t^8+2t^9+t^{13}$\\
\hline
$13$ & $1+t^2+2t^4+4t^6+4t^8+2t^{10}+t^{12}+t^{14}$\\
\hline
$14$ & $1+t^3+2t^4+4t^7+4t^8+2t^{11}+t^{12}+t^{15}$\\
\hline
$15$ & $1+t^2+3t^4+6t^6+7t^8+6t^{10}+3t^{12}+t^{14}+t^{16}$\\
\hline
$16$ & $1+2t^4+3t^5+8t^8+8t^9+3t^{12}+2t^{13}+t^{17}$\\
\hline
$17$ & $1+t^2+3t^4+8t^6+12t^8+12t^{10}+8t^{12}+3t^{14}+t^{16}+t^{18}$\\
\hline
$18$ & $1+t^3+3t^4+10t^7+13t^8+5t^9+5t^{10}+13t^{11}+10t^{12}+3t^{15}+t^{16}+t^{19}$\\
\hline
$19$ & $1+t^2+3t^4+10t^6+20t^8+24t^{10}+20t^{12}+10t^{14}+3t^{16}+t^{18}+t^{20}$\\
\hline
$20$ & $1+3t^4+4t^5+t^6+4t^7+22t^8+28t^9+13t^{10}+13t^{11}+28t^{12}+22t^{13}+4t^{14}+$\\
& $t^{15}+4t^{16}+3t^{17}+t^{21}$\\
\hline
\end{tabular}
\end{table}

\begin{table}[ht] 
\caption{Hilbert series for $n=2$}
\begin{tabular} {|c|l|l|}
\hline
$k$& $H(\Lambda(S^kV)^{\OO(2)}, t)$ & $H(\Lambda(S^kV)^{\SO(2)}, t)$\\
\hline
$3$ & $1+ t^4$ & $1+ 2t^2 + t^4$\\
\hline
$4$ & $1+t+t^4+t^5$ & $1+t+2t^2+2t^3+t^4+t^5$\\
\hline
$5$ & $1+3t^4$ & $1+3t^2+3t^4+t^6$\\
\hline
$6$ & $1+t+t^3+4t^4+3t^5$ & $1+t+3 t^2+5t^3+5t^4+3t^5+t^6+t^7$\\
\hline
$7$ & $1+7t^4+t^8$ & $1+4 t^2+8t^4+4t^6+t^8$\\
\hline
$8$ & $1+t+2t^3+9t^4+9t^5+2t^6+$ & $1+t+4t^2+8t^3+12t^4+12t^5+8t^6+$\\
&$t^8+t^9$ &$4t^7+t^8+t^9$\\
\hline
$9$ & $1+14t^4+4t^6+5t^8$ & $1+5t^2+18t^4+18t^6+5t^8+t^{10}$\\
\hline
$10$ & $1+t+4t^3+17t^4+21t^5+11t^6+$ & $1+t+5t^2+13t^3+24t^4+32t^5+32t^6+$\\
& $7t^7+9t^8+5t^9$ & $24t^7+13t^8+5t^9+t^{10}+t^{11}$\\
\hline
$11$ & $1+24t^4+19t^6+24t^8+t^{12}$ & $1+6t^2+33t^4+58t^6+33t^8+6t^{10}+t^{12}$\\
\hline
$12$ & $1+t+6t^3+29t^4+44t^5+37t^6+$ & $1+t+6t^2+18t^3+43t^4+73t^5+94t^6+$\\
& $37t^7+44t^8+29t^9+6t^{10}+$ & $94t^7+73t^8+43t^9+18t^{10}+6t^{11}+$\\
&$t^{12}+t^{13}$ &$t^{12}+t^{13}$\\
\hline
$13$ & $1+38t^4+58t^6+93t^8+$ & $1+7t^2+55t^4+151t^6+151t^8+55t^{10}+$\\
&$17t^{10}+7t^{12}$ &$7t^{12}+t^{14}$\\
\hline
$14$ & $1+t+9t^3+45t^4+81t^5+96t^6+$ & $1+t+7t^2+25t^3+69t^4+141t^5+227t^6+$\\
& $127t^7+162t^8+131t^9+60t^{10}+$ & $289t^7+289t^8+227t^9+141t^{10}+69t^{11}+$\\
&$24t^{11}+16t^{12}+7t^{13}$ & $25t^{12}+7t^{13}+t^{14}+t^{15}$\\
\hline
\end{tabular}
\end{table}


\begin{table}[ht] \label{table3}
\caption{Hilbert series for $n=3$}
\begin{tabular} {|c|l|}
\hline
$k$& $H(\Lambda(S^kV)^{\SL(3)}, t)$ \\
\hline
$3$ & $1+t^3+t^7+t^{10}$\\
\hline
$4$ & $1+t^6+t^9+t^{15}$\\
\hline
$5$ & $1+t^3+t^6+t^9+t^{12}+t^{15}+t^{18}+t^{21}$\\
\hline
$6$ & $1+t^5+2t^6+t^7+t^8+6t^9+7t^{10}+6t^{11}+8t^{12}+13t^{13}+16t^{14}+13t^{15}+8t^{16}+$\\
& $6t^{17}+7t^{18}+6t^{19}+t^{20}+t^{21}+2t^{22}+t^{23}+t^{28}$\\
\hline
\end{tabular}
\end{table}

\begin{table}[ht] 
\caption{Hilbert series for $n=3$}
\begin{tabular} {|c|l|l|}
\hline
$k$& $H(\Lambda(S^kV)^{\OO(3)}, t)$ & $H(\Lambda(S^kV)^{\SO(3)}, t)$\\
\hline
$3$ & $1+2t^4+2t^6+t^{10}$ & $1+3t^3+2t^4+2t^6+3t^7+t^{10}$\\
\hline
$4$ & $1+t+3t^4+12t^5+15t^6+8t^7+8t^8+$ & $1+t+3t^4+12t^5+15t^6+8t^7+8t^8+$\\
& $15t^9+12t^{10}+3t^{11}+t^{14}+t^{15}$ & $15t^9+12t^{10}+3t^{11}+t^{14}+t^{15}$\\
\hline
$5$ &$1+10t^4+60t^6+158t^8+294t^{10}+$ & $1+7t^3+10t^4+15t^5+60t^6+125t^7+$ \\
& $210t^{12}+125t^{14}+15t^{16}+7t^{18}$ & $158t^8+210t^9+294t^{10}+294t^{11}+$\\
& &$210t^{12}+158t^{13}+125t^{14}+60t^{15}+$\\
& &$15t^{16}+10t^{17}+7t^{18}+t^{21}$\\
\hline
\end{tabular}
\end{table}

\begin{table}[ht] 
\caption{Hilbert series for $n=4$}
\begin{tabular} {|c|l|l|}
\hline
$k$& $H(\Lambda(S^kV)^{\SL(4)}, t)$ & $H(\Lambda(S^kV)^{\Sp(4)}, t)$\\
\hline
$3$ & $1+t^4+t^8+t^{12}+t^{16}+t^{20}$ & $1+t^2+3t^4+4t^6+7t^8+7t^{10}+$\\
&& $7t^{12}+4t^{14}+3t^{16}+t^{18}+t^{20}$\\
\hline
$4$ & $1+t^5+t^8+t^9+t^{10}+4t^{13}+$ & $1+t^4+5t^5+2t^6+2t^7+22t^8+55t^9+$\\
& $3t^{14}+t^{15}+3t^{16}+4t^{17}+4t^{18}+$ & $79t^{10}+106t^{11}+207t^{12}+383t^{13}+$\\
& $3t^{19}+t^{20}+3t^{21}+4t^{22}+t^{25}+$ & $490t^{14}+528t^{15}+651t^{16}+804t^{17}+$\\
& $t^{26}+t^{27}+t^{30}+t^{35}$& $804t^{18}+651t^{19}+528t^{20}+490t^{21}+$\\
& & $383t^{22}+207t^{23}+106t^{24}+79t^{25}+55t^{26}+$\\
& & $22t^{27}+2t^{28}+2t^{29}+5t^{30}+t^{31}+t^{35}$\\
\hline
\end{tabular}
\end{table}

\begin{table}[ht] 
\caption{Hilbert series for $n=4$}
\begin{tabular} {|c|l|l|}
\hline
$k$& $H(\Lambda(S^kV)^{\OO(4)}, t)$ & $H(\Lambda(S^kV)^{\SO(4)}, t)$\\
\hline
$3$ & $1+2t^4+3t^6+15t^8+8t^{10}+15t^{12}+$ & $1+7t^4+4t^6+30t^8+16t^{10}+30t^{12}+$\\
& $t^{14}+5t^{16}$ & $4t^{14}+7t^{16}+t^{20}$\\
\hline
\end{tabular}
\end{table}

\begin{example} \label{ex_Hilb5} Hilbert series of $\Lambda(S^3V)^G$ for $n=5$:
\[
H(\Lambda(S^3V)^{\SL(5)}, t) = 1+t^5+t^{30}+t^{35};
\]
\begin{align*}
H(\Lambda(S^3V)^{\OO(5)}, t)= &1+2t^4+3t^6+24t^8+79t^{10}+198t^{12}+461t^{14}+592t^{16}+ \\
&741t^{18}+488t^{20}+359t^{22}+106t^{24}+56t^{26}+3t^{28}+7t^{30};
\end{align*}
\begin{align*}
&H(\Lambda(S^3V)^{\SO(5)}, t) = 1+2t^4+7t^5+3t^6+3t^7+24t^8+56t^9+79t^{10}+106t^{11}+\\
&198t^{12}+359t^{13}+461t^{14}+488t^{15}+592t^{16}+741t^{17}+741t^{18}+592t^{19}+\\
&488t^{20}+461t^{21}+359t^{22}+198t^{23}+106t^{24}+79t^{25}+56t^{26}+24t^{27}+3t^{28}+\\
&3t^{29}+7t^{30}+2t^{31}+t^{35}.
\end{align*}
\end{example}

\begin{example}\label{ex_LambdaS3_SL_inv} 
The results in Table 3 show that $\Lambda(S^3V)^{\SL(3)}$ is generated by a pair $\{v, \ast{v}\}$, where $v$ is a certain element in $\Lambda^3(S^3V)$ and $\ast{v}$ is the Hodge dual of $v$, i.e. the unique element in $\Lambda^{7}(S^3V)$ such that
$$
v \wedge \ast{v} = \left\langle v, v \right\rangle \mathrm{vol},
$$
where $\left\langle \cdot, \cdot \right\rangle$ is the inner product on $\Lambda^3(S^3V)$, which is induced by the standard inner product on $V = \CC^3$, and $\mathrm{vol}$ is a chosen volume form on $S^3V$. In this example, we give explicit expressions for $v$ and $\ast{v}$.

Let $\{e_1, e_2, e_3\}$ denote the standard basis for $V = \CC^3$. Then a basis for $V^{\otimes 3}$ is given by the elements $\{e_i\otimes e_j\otimes e_k\}_{\{1\leq i,j,k \leq 3\}}$. We consider the standard symmetrization mapping on $V^{\otimes 3}$, defined by
\begin{align*}
&S: V^{\otimes 3} \rightarrow V^{\otimes 3}\\
&S(e_i \otimes e_j \otimes e_k) = \frac{1}{3!} \sum_{\sigma \in S_3} e_{\sigma(i)}\otimes e_{\sigma(j)} \otimes e_{\sigma(k)},
\end{align*}
where $S_3$ denotes again the symmetric group in $3$ variables. Using the standard notation
$e_ie_je_k = S(e_i \otimes e_j \otimes e_k)$, a basis for $S^3V = \mathrm{im}S$ is given by the following ten elements:
\begin{align*}
&a_1 = e_1^3, \quad a_2 = e_2^3, \quad a_3 = e_3^3, \quad a_4 = e_1^2e_2, \quad a_5 = e_1^2e_3, \\
&a_6 =e_2^2e_3, \quad a_7= e_1e_2^2, \quad a_8 = e_1e_3^2, \quad a_9 = e_2e_3^2, \quad a_{10} = e_1e_2e_3.
\end{align*}

Next, the following map 
$$
f: V^{\otimes 3p} \rightarrow \Lambda^p(S^3V)
$$ 
given by
$$
e_{i_1}\otimes e_{j_1} \otimes e_{k_1} \otimes \cdots \otimes e_{i_p}\otimes e_{j_p} \otimes e_{k_p} \mapsto e_{i_1}e_{j_1}e_{k_1} \wedge \cdots \wedge e_{i_p}e_{j_p}e_{k_p}
$$
is a surjective homomorphism of $\GL(3)$-modules for any $p \geq 1$. Thus, any $\SL(3)$-invariant in $\Lambda^p(S^3V)$ comes from an $\SL(3)$-invariant in $V^{\otimes 3p}$. In particular, for $p=3$, we take $v$ to be the image of the element $\mathrm{Std_3}(e_1, e_2, e_3)\otimes\mathrm{Std_3}(e_1, e_2, e_3)\otimes \mathrm{Std_3}(e_1, e_2, e_3)$, where
$$
\mathrm{Std_3}(e_1, e_2, e_3) = \sum_{\sigma \in S_3} \mathrm{sign}(\sigma) e_{\sigma(1)}\otimes e_{\sigma(2)}\otimes e_{\sigma(3)}
$$
is the standard polynomial of degree $3$. Explicitly, we obtain
\begin{align*}
v = &a_1 \wedge a_2 \wedge a_3 - 3a_3 \wedge a_4 \wedge a_7 - 3 a_1 \wedge a_6 \wedge a_9  + 3a_2 \wedge a_5 \wedge a_8 +\\
 &6a_7 \wedge a_8 \wedge a_{10} - 6a_4 \wedge a_9 \wedge a_{10} + 6 a_5 \wedge a_6 \wedge a_{10} +\\
&3a_5 \wedge a_7 \wedge a_9 + 3a_4 \wedge a_6 \wedge a_8.
\end{align*}
Then, if we choose $\mathrm{vol} = a_1 \wedge a_2 \wedge \cdots \wedge a_{10}$, for $\ast{v}$ we obtain
\begin{align*}
\ast{v} = &a_4 \wedge a_5 \wedge a_6 \wedge a_7 \wedge a_8 \wedge a_9 \wedge a_{10} - \frac{1}{3}a_1 \wedge a_2 \wedge a_5 \wedge a_6 \wedge a_8 \wedge a_9 \wedge a_{10} -\\
&-\frac{1}{3}a_2 \wedge a_3 \wedge a_4 \wedge a_5 \wedge a_7 \wedge a_8 \wedge a_{10} - \frac{1}{3}a_1 \wedge a_3 \wedge a_4 \wedge a_6 \wedge a_7 \wedge a_9 \wedge a_{10} -\\ 
&-\frac{1}{9}a_1 \wedge a_2 \wedge a_3 \wedge a_4 \wedge a_5 \wedge a_6 \wedge a_{9} + \frac{1}{9}a_1 \wedge a_2 \wedge a_3 \wedge a_5 \wedge a_6 \wedge a_7 \wedge a_{8} -\\
&-\frac{1}{9}a_1 \wedge a_2 \wedge a_3 \wedge a_4 \wedge a_7 \wedge a_8 \wedge a_{9} - \frac{1}{9}a_1 \wedge a_2 \wedge a_3 \wedge a_4 \wedge a_6 \wedge a_8 \wedge a_{10} +\\
&+\frac{1}{9}a_1 \wedge a_2 \wedge a_3 \wedge a_5 \wedge a_7 \wedge a_9 \wedge a_{10}.
\end{align*}
\end{example}

\begin{remark} Similarly, the results in Tables 1, 3, and 5 show that $\Lambda(S^5V)^{\SL(2)}$, $\Lambda(S^6V)^{\SL(2)}$, $\Lambda(S^8V)^{\SL(2)}$, $\Lambda(S^4V)^{\SL(3)}$, and $\Lambda(S^3V)^{\SL(5)}$ are generated each by a pair $\{v, \ast{v}\}$ and these pairs are obtained in an analogous way to Example \ref{ex_LambdaS3_SL_inv}.
\end{remark}

In the next set of examples, we consider $W = \Lambda^3V$. Then, Equation (\ref{eq_HilbSeriesLambda2}) implies that
\[
H(\Lambda (\Lambda^3V), x_1, \dots, x_n,t) = \prod_{1\leq i < j < j \leq n}(1+x_ix_jx_kt).
\]

In Tables 7-9 below, we give the expressions for $H(\Lambda (\Lambda^3V)^{G}, t)$, for $n=5,6$.

\begin{table}[ht] 
\caption{Hilbert series for $n=5$}
\begin{tabular} {|l|l|l|}
\hline
$H(\Lambda (\Lambda^3V)^{\SL(5)}, t)$ & $H(\Lambda (\Lambda^3V)^{\OO(5)}, t)$ & $H(\Lambda (\Lambda^3V)^{\SO(5)}, t)$ \\
\hline
$1+t^{10}$ &$1+t^{10}$& $1+t^3+t^7+t^{10}$\\
\hline
\end{tabular}
\end{table}

\begin{table}[ht] 
\caption{Hilbert series for $n=6$}
\begin{tabular} {|l|l|}
\hline
$H(\Lambda (\Lambda^3V)^{\SL(6)}, t)$ & $H(\Lambda (\Lambda^3V)^{\Sp(6)}, t)$ \\
\hline
$1+t^2+t^4+t^6+t^8+t^{10}+t^{12}+$& $1+2t^2+4t^4+5t^6+5t^8+5t^{10}+5t^{12}+$\\
$t^{14}+t^{16}+t^{18}+t^{20}$ & $5t^{14}+4t^{16}+2t^{18}+t^{20}$\\
\hline
\end{tabular}
\end{table}

\begin{table}[ht] 
\caption{Hilbert series for $n=6$}
\begin{tabular} {|l|l|}
\hline
$H(\Lambda (\Lambda^3V)^{\OO(6)}, t)$ & $H(\Lambda (\Lambda^3V)^{\SO(6)}, t)$ \\
\hline
$1+t^4+2t^8+3t^{10}+2t^{12}+t^{16}+t^{20}$ & $1+t^2+t^4+2t^6+2t^8+8t^{10}+2t^{12}+$\\
& $2t^{14}+t^{16}+t^{18}+t^{20}$\\
\hline
\end{tabular}
\end{table}

\begin{remark} In a similar way, further results for $H(\Lambda({S^kV})^G, t)$ and $H(\Lambda({\Lambda^kV})^G, t)$ for larger $k$ and for larger $n$ can be obtained. These results are not included in the paper only in order to avoid too complicated expressions for the respective Hilbert series.
\end{remark}

\section{Appendix}
 Below we provide the source code of the program in Mathematica for determining the Hilbert series $H(\Lambda(S^kV)^G, t)$ for $V = \CC^2$ and $G = \SL(2)$, $\OO(2)$, and $\SO(2)$. The value of $k$ should be chosen in the beginning of the program.

\begin{verbatim}
k = 3; H = x - y; 
For[i = 0, i <= k, i++, H = Expand[H*(1 + x^(k - i)*y^i*t)]];  
M = CoefficientRules[H, {x, y}]; L = {}; 
For[i = 1, i <= Length[M], i++, 
 If[M[[i]][[1]][[1]] > M[[i]][[1]][[2]], L = Append[L, M[[i]]]]]; 
MultSeries = FromCoefficientRules[L, {x, y}]; 
MultSeries = Expand[MultSeries*1/x]; 
Print["Multiplicity Series = " MultSeries];
M1 = CoefficientRules[MultSeries, {x, y}]; L1 = {}; 
For[i = 1, i <= Length[M1], i++, 
 If[M1[[i]][[1]][[1]] == M1[[i]][[1]][[2]], L1 = Append[L1, M1[[i]]]]]; 
InvariantsSL = FromCoefficientRules[L1, {x, y}]; 
InvariantsSL = InvariantsSL /. x -> 1; 
InvariantsSL = InvariantsSL /. y -> 1; 
Print["Hilbert series of SL-invariants =" InvariantsSL];
L2 = {}; For[i = 1, i <= Length[M1], i++, 
 If[EvenQ[M1[[i]][[1]][[1]]] && EvenQ[M1[[i]][[1]][[2]]], 
  L2 = Append[L2, M1[[i]]]]]; 
InvariantsO = FromCoefficientRules[L2, {x, y}]; 
InvariantsO = InvariantsO /. x -> 1; 
InvariantsO = InvariantsO /. y -> 1; 
Print["Hilbert series of O-invariants =" InvariantsO];
For[i = 1, i <= Length[M1], i++, 
  If[OddQ[M1[[i]][[1]][[1]]] && OddQ[M1[[i]][[1]][[2]]], 
   L2 = Append[L2, M1[[i]]]]]; 
InvariantsSO = FromCoefficientRules[L2, {x, y}]; 
InvariantsSO = InvariantsSO /. x -> 1; 
InvariantsSO = InvariantsSO /. y -> 1; 
Print["Hilbert series of SO-invariants =" InvariantsSO]
\end{verbatim}


\section*{Acknowledgements}
I am very grateful to Vesselin Drensky for introducing me to the topic of Hilbert series and for his advices and stimulating discussions.

\end{document}